\documentclass[12pt]{amsart}
\usepackage{amsfonts,amsmath,oldgerm,amssymb,amscd}
\newtheorem{theorem}{Theorem}
\newtheorem{lemma}[theorem]{Lemma}
\newtheorem{corollary}{Corollary}

\theoremstyle{definition}                 
\newtheorem*{example}{Example} 

\theoremstyle{definition}

\usepackage{amsfonts}
\newcommand{\field}[1]{\mathbb{#1}}          \newcommand{\Q}{\field{Q}}
                   \newcommand{\Z}{\field{Z}}
\newcommand{\C}{\field{C}}

\newcommand{\ord}{{\rm ord}}
\newcommand{\mO}{{\mathcal O}}
\newcommand{\mf}{\mathfrak}
\newcommand{\us}{\underset}
\newcommand{\os}{\overset}
\newcommand{\card}{{\rm card}}
\newcommand{\rank}{{\rm rank}}

\newcommand{\cgtu}{\hskip-13pt^u\hskip+10pt}

\begin{document}

\title[${\bf S}$-Arithmetic Groups of ${\bf SL_2}$ type]
{${\bf S}$-Arithmetic Groups of ${\bf SL_2}$ type}

\author{R. Sarma}

\address{Harish-Chandra Research Institute, Jhunsi, Chhatnag road, Allahabad,
India.}

\email{ritumoni@mri.ernet.in}

\subjclass{Primary  20F05,  11F06;  Secondary  22E40\\{\bf  R. Sarma.}
Harish-Chandra Research Institute,  Chhatnag Road, Jhunsi,  Allahabad  211 019,
India. ritumoni@mri.ernet.in}
\date{}

\begin{abstract}
For a number field $K$, we show that any $S$-arithmetic subgroup
of $SL_2(K)$ contains a subgroup of finite index generated by three 
elements if $\card(S)\ge 2$.
\end{abstract}

\maketitle

\section{Introduction and Notation}

Let $K$ be a number field and let $S_\infty$ be the set of all 
nonconjugate embeddings of $K$ into $\C$. We refer to these embeddings as
{\it infinite primes} of $K$. 
If $r_1$ (resp. $r_2$) is the number of distinct 
real (resp. nonconjugate complex) embeddings so that the cardinality of 
$S_\infty$ is $r_1+r_2$, then $r_1+2r_2=[K:\Q]$, the extension degree of
$K$. Let $S$ be a finite set of primes in $K$
containing $S_\infty$.  The ring of integers in $K$ is denoted by $\mO_K$.
For a prime ideal $\mf{p}$ of $\mO_K$, 
denote by $v_{\mf{p}}$ the valuation
defined by $\mf{p}$.
The ring $\mO_S:=\{x\in K: v_{\mf{p}}(x)\ge 0 ~~{\rm for~~ every~~
prime }~~\mf{p}\notin S\}$ is called the ring of {\it $S$-integers} of $K$. 
Then $\mO_{S_\infty}=\mO_K$.
If $F$ is a subfield of $K$, then set
\begin{equation}\label{definition1}
S(F):=\{\mf{p}\cap\mO_F:\mf{p}\in S-S_\infty\}\sqcup S_\infty(F)
\end{equation}
where $S_\infty(F)$ denotes the infinite primes of $F$.
We write 
\begin{equation}\label{definition2}
\mO_{S(F)}:=\{x\in F: v_\mf{p}(x)\ge0~\forall\,\mf{p}\notin S(F)\}
\end{equation}
the ring of $S(F)$-integers in $F$.

For two subgroups $H_1$ and $H_2$ in a group,
if $H_1\cap H_2$ is a subgroup of finite index both in $H_1$ and $H_2$, then
we say that $H_1$ and $H_2$ are {\it commensurable} and we write 
$H_1\asymp H_2$. In particular, a group is
commensurable with its subgroups of finite index. 
A subgroup $\Gamma$ of $SL_2(K)$ is called {\it $S$-arithmetic}
if $\Gamma\asymp SL_2(\mO_S)$.


A subset $X$ of a group $G$ is called a set of {\it virtual} generators
of $G$ if the group generated by $X$ is a subgroup of finite index in $G$ and 
the group $G$ is said to be generated {\it virtually} by $X$. 

Let the cardinality of any set $X$ be denoted by $\card(X)$. 

A number field is called a {\it totally real} field if all its embeddings
are real. A number field is called a {\it CM field} if it is an imaginary 
quadratic extension of a totally real field. If a number field is 
not CM then we refer to it as a {\it non-CM} field.

For any commutative ring $A$, denote by 
\begin{equation}\begin{pmatrix}1& A\cr0&1\end{pmatrix}
\phantom{....}(\textnormal{ resp. }\begin{pmatrix}1&0\cr A&1\end{pmatrix})
\end{equation}
the subgroup of $SL_2(A)$ consisting
of matrices of the form 
\begin{equation*}
\begin{pmatrix}1&x\cr0&1\end{pmatrix} 
\phantom{....}(\textnormal{ resp. }\begin{pmatrix}1&0\cr x&1\end{pmatrix}) 
\phantom{......}\textnormal{ for } x\in A.
\end{equation*}

Let $G$ be any group and let $a,b\in G$. 
Denote by $^ab$ the element $aba^{-1}$ in $G$.


We use a few well known number theoretic results 
(for details, see \cite{Marcus},\cite{Platonov}): 
The ring $\mO_K$ of integers in $K$
is a {\bf Dedekind domain}. An ideal of $\mO_K$ has a unique 
{\bf factorization into prime ideals} of $\mO_K$. 
For a finitely generated abelian group $H$, let $\rank(H)$ denote the 
{\bf rank} of $H$ as a $\Z$-module. {\bf Dirichlet's unit theorem} asserts that
\begin{equation}\label{rankformula1}
\rank(\mO_K^*)=r_1+r_2-1
\end{equation}
where $r_1$ and $r_2$ are defined
as above. Also  (cf. Lemma \ref{OS})
\begin{equation}\label{rankformula2}
\rank(\mO_S^*)=\card(S)-1.
\end{equation}

The group of units of a ring $A$ is denoted by $A^*$. For an ideal
$\mf{a}$ of $\mO_K$, let the {\bf order} of the class of $\mf{a}$ in the 
{\bf ideal class group} of $K$ be denoted by $\ord_K(\mf{a})$. 
It is well known that the class group of a number field is finite. Thus 
$\ord_K(\mf{a})$ is always a {\bf finite} number.

Now we state the main result of the paper.
\begin{theorem}\label{maintheorem}
Let $K$ be a number field and let $S$ be a finite set 
of primes in $K$ containing the infinite ones such that ${\rm card}(S)\ge2$. 
Any $S$-arithmetic subgroup of $SL_2(K)$ is virtually generated by three 
elements.
\end{theorem}
We postpone the proof of this theorem to section 3. 
In \cite{mine2}, it is shown that the {\bf higher rank arithmetic groups} are 
virtually generated by three elements. The tools used to prove this do not
seem to work for the case of $S$-arithmetic groups. For instance,
if $U$ is a {\bf unipotent group}, and if $\Gamma$ is 
a {\bf Zariski dense} subgroup of an
arithmetic subgroup of $U$, then $\Gamma$ 
is also arithmetic. This fact plays a crucial role in the case of
higher rank arithmetic groups. The analogous statement does not hold in the 
case of $S$-arithmetic subgroups. So it needs a separate
treatment. The case of $SL_2$ is the first case that one would like to deal 
with because this is the simplest possible case. The technique here may 
indicate 
how to proceed for other $S$-arithmetic groups. In the next section we prove
a number theoretic result which asserts that $\mO_S$ is {\bf almost} generated
by some unit (in fact, by any positive power of that unit) in $\mO_S$.  
Then our main result follows from a theorem due to Vaserstein.
The condition that $\card(S)\ge2$ is equivalent to say that the group 
$\mO_S^*$ is infinite. 


\section{A Number Theoretic Result}\label{section2}
\begin{theorem}\label{goodunit}
Let $K$ be a non-CM field and let $S$ be a finite set of primes including the
infinite ones with ${\rm card}(S)\ge 2$. Then there exists $\alpha\in \mO_S^*$ 
such that the ring $\Z[\alpha^n]$ is a subgroup of finite index in the
ring $\mO_S$ of $S$-integers for every positive integer $n$.
\end{theorem}

The proof of Theorem \ref{goodunit} is divided into a few lemmata. 
\begin{lemma} \label{infiniteindex}
If $K$ is a non-CM field and if $F$ is a proper subfield of $K$, then 
$\mO_F^*$ is a subgroup of infinite index in 
$\mO_K^*$.
\end{lemma}

\noindent
\begin{proof}
Assume to the contrary that there exists a subfield $F$ of $K$ such that 
the quotient $\mO_K^*/\mO_F^*$ is finite.
We show that $K$ is a CM field and arrive at a contradiction.

Let $d=[K:F]$ be the degree of exention. 
Let $A$ (resp. $B$) be the set of all real (resp. nonconjugate complex)
embeddings of $F$ over $\Q$.
For an embedding $a\in A$ \index{embedding!real}\index{embedding!complex},
let $x(a)$ (resp. $y(a)$) be the number of real (resp. nonconjugate complex)
extensions of $a$ to $K$. Then we have 
\begin{equation}\label{realcomplexembeddings}
x(a)+2y(a)=d. 
\end{equation}
Similarly, we define 
$x(b),y(b)$ for $b\in B$. Then $x(b)=0$ and $y(b)=d$. Since
$\mO_K^*/\mO_F^*$ is finite, $\rank(\mO_F^*)=\rank(\mO_K^*)$. Hence by
\eqref{rankformula1}, we have
\begin{eqnarray} \label{cardinality}
\card(A)+\card(B)-1 & = &\us{a\in A}{\sum}x(a)+\big(\us{a\in A}{\sum}y(a)+
\us{b\in B}{\sum}y(b)\big)-1 \nonumber \\
& = &\us{a\in A}{\sum}\{x(a)+y(a)\}+\us{b\in B}{\sum}y(b)-1 \nonumber\\
& \ge &\card(A)+d~\card(B)-1, 
\end{eqnarray}
as $x(a)+y(a)\ge 1$. Now since $F$ is a proper subfield of $K$, we have
$d> 1$. Hence using \eqref{cardinality}, we see that $\card(B)= 0$. 
Thus, $F$ is a totally real field. 
Therefore by inequality \eqref{cardinality}, we get 
$\card(A)= \us{a\in A}{\sum}\{x(a)+y(a)\}$. Hence for each $a\in A$, 
we have 
\begin{equation}\label{realextensions}
x(a)+y(a)=1.
\end{equation}
Now, it follows from \eqref{realcomplexembeddings} that
$y(a)=d-1$. Therefore, since $d-1\ge1$, by \eqref{realextensions},
we have $x(a)=0$ and $y(a)=1$.
Thus again by \eqref {realcomplexembeddings}, we get $d=2$ so that
the field $K$ is an imaginary quadratic extension over the totally real field
$F$ which is a contradiction as desired.
\end{proof}

\begin{lemma}\label{firstcomparison}
Let $K=\Q(\alpha)$ and let $\alpha$ be integral. Then $\Z[\alpha^{-1}]$
is of finite index in $\mO_K[\alpha^{-1}]$.
\end{lemma}
\begin{proof}
Since $\alpha$ is an integral element, we have 
$\Z[\alpha]\subset \Z[\alpha^{-1}]$.
Let $n$ be the index of $\alpha\mO_K$ in $\mO_K$. For $0\le i\le (n-1)$,
the cosets $\alpha\mO_K+i$ are the distinct cosets. 
Otherwise, $\alpha\mO_K+i=\alpha\mO_K+j$ for $0\le i<j\le(n-1)$ 
so that $j-i\in\alpha\mO_K$. This implies that $n$ divides $j-i$
which is a contradiction. Therefore, $\mO_K$ is the union of these $n$
cosets. Thus, in particular, 
\begin{equation}\label{twist}
\Z[\alpha]+\alpha\mO_K=\mO_K.
\end{equation}
On the other hand, $\Z[\alpha]$ is of finite index in $\mO_K$. 
By \eqref{twist}, we may assume that the distinct cosets 
(as an additive subgroup) of $\Z[\alpha]$
in $\mO_K$ are $\Z[\alpha]+\alpha x_i$ with $0\le i\le (m-1)$. 
We claim that the representatives of $\mO_K[\alpha^{-1}]/\Z[\alpha^{-1}]$
in $\mO_K[\alpha^{-1}]$, 
are $\alpha x_i$ (not necessarily  distinct). 
Let $y\in\mO_K$. Then, by \eqref{twist}, 
$y=y_1+\alpha x_{i_1}$ for $y_1\in\Z[\alpha]$ and $0\le i_1\le(m-1)$.
Thus $\alpha^{-1}y=\alpha^{-1}y_1+x_{i_1}$. Again, using \eqref{twist}, 
we have $x_{i_1}=z_1+\alpha x_{i_2}$ for $z_1\in\Z[\alpha]$ and 
$0\le i_2\le (m-1)$ so that 
$\alpha^{-1}y=(\alpha^{-1}y_1+z_1)+\alpha x_{i_2}$. Therefore,
$\Z[\alpha^{-1}]+\alpha^{-1}y=\Z[\alpha^{-1}]+\alpha x_{i_2}$. 
Thus inductively, one can show that 
$\Z[\alpha^{-1}]+\alpha^{-r}y=\Z[\alpha^{-1}]+\alpha x_i$ for some 
$0\le i\le(m-1)$.
\end{proof}

\begin{lemma}\label{OS}
Let $K$ be a number field and let $S$ be a finite set of primes in $K$ 
containing $S_\infty$.
Assume that
$S-S_\infty=\{\mf{q}_1,\dots,\mf{q}_r\}$, $\ord_K(\mf{q}_i)=a_i$ and that 
$\mf{q}_i^{a_i}$ is generated by $\beta_i\in\mO_K ~\forall\,i$.
Then $\mO_S=\mO_K[\beta_1^{-1},\dots,\beta_r^{-1}].$
\end{lemma}
\begin{proof}
Obviously, $\mO_S\supset\mO_K[\beta_1^{-1},\dots,\beta_r^{-1}].$
To see the other containment, let $x\in \mO_S$. Then $x=yz^{-1}$ for 
$y,z\in\mO_K$ and $v_\mf{p}(z)=0$ for $\mf{p}\notin S$ so that, 
by prime factorization,
$z\mO_K=\underset{i=1}{\overset{r}{\prod}}\mf{q}_i^{n_i}$ with $n_i\ge 0$.
Let $m=\underset{i=1}{\overset{r}{\prod}}a_i$. Since $\mf{q}_i^{a_i}$ 
is generated by $\beta_i$, we have
$z^{-m}=u\underset{i=1}{\overset{r}{\prod}}\beta_i^{-n_i'}$ for some 
$u\in\mO_K^*$ and $n_i'\ge0$ so that
$z^{-m}\in \mO_K[\beta_1^{-1},\dots,\beta_r^{-1}]$. 
Further, $z^{-1}=z^{m-1}z^{-m}$ and $z^{m-1}\in\mO_K$. Therefore,
$z^{-1}\in\mO_K[\beta_1^{-1},\dots,\beta_r^{-1}]$ and hence
$x=yz^{-1}\in\mO_K[\beta_1^{-1},\dots,\beta_r^{-1}]$. 
\end{proof}
Now by Lemma \ref{firstcomparison} and Lemma \ref{OS}, 
we have the following lemma.
\begin{lemma}\label{comparison}
Suppose that
$R$ is a subring of finite index in $\mO_K$. Then with the notation as in 
Lemma \ref{OS}, the ring $R[\beta_1^{-1},\dots,\beta_r^{-1}]$ is of finite 
index in $\mO_S$.
\hfill$\square$
\end{lemma}

\noindent
Let $\{S_i:1\le i\le s\}$ be
the set of {\bf all the proper subsets} of $S$ and let 
$\{K_j:1\le j\le t\}$ be 
the set of {\bf all the proper subfields} of $K$. Define 
\begin{equation}\label{definition3}
V_i:=\mO_{S_i}^*\otimes_{\Z}\Q~~~~\textnormal{ and }~~~~
W_j:=(\mO_{S(K_j)}^*\cap\mO_S^*)\otimes_{\Z}\Q.
\end{equation}
Then $V_i$ (resp. $W_j$) is vector  space of dimension $\rank(\mO_{S_i}^*)$
(resp. $\rank(\mO_{S(K_j)}^*)$) over $\Q$.
By Lemma \ref{OS}, we have
$\mO_S^*=\mO_K^*\times \Z^r$ where $r=\card(S)-\card(S_\infty)$.

\begin{lemma}\label{step1}
With the above notation, if $K$ is a non-CM field, there exists 
$\alpha\in\mO_S^*-(\os{s}{\us{i=1}{\cup}}V_i)
\cup(\os{t}{\us{j=1}{\cup}}W_j)$
such that $v_\mf{p}(\alpha)<0$ for all $\mf{p}\in S-S_\infty$.
\end{lemma}

\begin{proof} 
For each $1\le j\le s$, we have 
(see \eqref{definition1} and \eqref{definition2} for definition)
\begin{eqnarray}
{\rm rank}(\mO_{S(K_j)}^*)&=&\card(S(K_j))-1\nonumber
\\&=& \{\card(S_\infty(K_j))-1\}+{\rm card}(S(K_j)-S_\infty(K_j))\nonumber\\
&=& \rank(\mO_{K_j}^*)+{\rm card}(S(K_j)-S_\infty(K_j)).
\end{eqnarray}
Since $K$ is a non-CM field, by Lemma \ref{infiniteindex}, 
$\rank(\mO_{K_j}^*)<\rank(\mO_K^*)$.
Moreover,
$\card(S(K_j)-S_\infty(K_j))\le\card(S-S_\infty)$. Therefore, we get
\begin{eqnarray}\label{inequality}
     {\rm rank}(\mO_{S(K_j)}^*\cap\mO_S^*)<{\rm rank}(\mO_S^*).
\end{eqnarray}
Further, ${\rm rank}(\mO_{S_i}^*)={\rm card}(S_i)-1<{\rm rank}(\mO_S^*)$.
Then by comparing the dimensions, we have $V_i\subsetneqq V$ and 
$W_j\subsetneqq V$ where $V:=\mO_S^*\otimes_\Z\Q$
(cf. \eqref{definition3}).
Since a finite union of proper subspaces of a vector space
over an infinite field is a proper subset of the vector space, 
we have  $V-(\os{s}{\us{i=1}{\cup}}V_i)
\cup(\os{t}{\us{j=1}{\cup}}W_j)$ 
is {\bf nonempty}. Let
\begin{equation}\label{definitionX}
X:=\{x\in \mO_S^*:\,v_\mf{p}(x)~<0\,\forall\mf{p}\in S-S_\infty\}.
\end{equation}
Then $X$ is Zariski dense in $V$. 
Thus the set $$Y:=X-(\os{s}{\us{i=1}{\cup}}V_i)
\cup(\os{t}{\us{j=1}{\cup}}W_j)$$ is also nonempty.
If $\alpha\in Y$, then $\alpha^n\in Y$.
Thus, $\alpha\in \mO_S^*$ can be chosen such that $v_\mf{p}(\alpha)<0$
for each $\mf{p}\in S-S_\infty$.
\end{proof}

\begin{lemma}Let $\alpha$ be chosen as in Lemma \ref{step1}. Then
the ring $\Z[\alpha^n]$ is a subgroup of finite index in $\mO_S$ for every
positive integer $n$.
\end{lemma}

\begin{proof}
We claim $\Q(\alpha)=K$. If not, then let $\Q(\alpha)=L$ such that
$L\subsetneqq K$. Assume 
for $\mf{p}\notin S$ and $x\in \mO_L$ that
$v_{\mf{p}\cap\mO_L}(x)\ne0$ so that 
$x\mO_L\subset\mf{p}\cap\mO_L$. Then, 
$x\mO_K\subset(\mf{p}\cap\mO_L)\mO_K\subset\mf{p}$ so that
$v_{\mf{p}}(x)\ne0$. Thus, equivalently, for $x\in\mO_L$,
if $v_{\mf p}(x)=0$ for every $\mf{p}\notin S$, we have $v_\mf{p}(x)=0$ 
for every $\mf{p}\notin
S(L)$. Therefore, in particular, $v_p(\alpha^{-1})=0\,\forall\mf{p}
\notin S(L)$
so that $\alpha\in \mO_{S(L)}^*\cap \mO_S^*$.
This contradicts the choice of $\alpha$. Hence $\Q(\alpha)=K$. 

Since $K=\Q(\alpha)$, we also have $K=\Q(\alpha^{-1})$ and since
$\alpha^{-1}$ is integral in $K$, the ring $\Z[\alpha^{-1}]$
is a subgroup of finite index in $\mO_K$. 
Let $S-S_\infty=\{\mf{p}_i:1\le i\le l\}$. 
Consider the prime factorization
\begin{equation}
\alpha^{-1}\mO_K=\underset{i=1}{\overset{l}{\prod}}\mf{p}_i^{n_i}
\end{equation}
where $n_i\ge0$. In fact, by our choice of $\alpha$, 
$n_i>0$ for each $i$. 
Let $\ord_K(\mf{p}_i)=r_i$ and let $\mf{p}_i^{r_i}=\beta_i\mO_K$ 
for $\beta_i\in \mO_K$.
Then, we have 
\begin{equation}\label{alpha}
\alpha^m=\underset{i=1}{\overset{l}{\prod}} \beta_i^{-b_i}  
\end{equation} 
(if necessary absorbing a unit with $\beta_1$)
for some integers $m>0$ and $b_i>0$. 
Since $\beta_i\in\mO_K$,
it follows by \eqref{alpha} that $\beta_i^{-1}\in \mO_K[\alpha]$.
Now by Lemma \ref{comparison}, the ring 
$\mO_K[\alpha]=\mO_S$. Thus, by Lemma \ref{firstcomparison}, the ring 
$\Z[\alpha]$ is of finite index in $\mO_S$.
\end{proof}
This completes the proof of Theorem \ref{goodunit}. 

In fact, we have proved more. 
\begin{corollary}\label{goodunit2}
Let $K$ be any finite extension of $\Q$ and let $S$ be as before.
If ${\rm rank}(\mO_{S(L)}^*\cap\mO_S^*)<{\rm rank}(\mO_{S}^*)$ {\bf for every
proper subfield} $L$ of $K$, then there exists $\alpha\in\mO_S^*$
such that the ring $\Z[\alpha^n]$ is a subgroup of finite index in
$\mO_S$ {\bf for every $n\ge1$}.\hfill$\square$
\end{corollary}

\noindent
The hypothesis of Corollary \ref{goodunit2} {\bf may hold sometimes even for a 
CM field}. Here we see two examples:
\begin{example}
(i) The field $K=\Q(\sqrt{-1})$ is a CM field and  
$\mO_K=\Z[\sqrt{-1}]$.
The prime ideal $2\Z$ of $\Q$ is totally ramified in $K$. In fact, 
$2\mO_K=\mf{p}^2$ where $\mf{p}=\big<1+\sqrt{-1}\big>$. Let 
$S-S_\infty=\{\mf{p}\}$.
For $K$, the set $S_\infty$ of infinite primes is singleton. Thus 
${\rm card}(S)=2$ and 
hence ${\rm rank}(\mO_S^*)=1$. Also, $\mO_{S(\Q)}=\Z[\frac{1}{2}]$ and so
${\rm rank}(\mO_{S(\Q)}^*\cap\mO_S^*)=1$
(observe that 
$\mO_S=\Z[\sqrt{-1}][\frac{1}{1+\sqrt{-1}}]$ includes $\mO_{S(\Q)}$).
This is an example which does not satisfy the hypothesis 
corollary \ref{goodunit2}.

\noindent
(ii) Next consider the ideal $5\Z$ of $\Q$ which splits completely in $K$:
$5\mO_K=\mf{p}_1\mf{p}_2$ where $\mf{p}_1=\big<5,2+\sqrt{-1}\big>$ and 
$\mf{p}_2=\big<5,2-\sqrt{-1}\big>$. Let $S-S_\infty=\{\mf{p}_1,\mf{p}_2\}$. 
Then ${\rm card}(S)=3$ and hence ${\rm rank}(\mO_S^*)=2$. The contraction of
the primes of $S-S_\infty$ to $\Q$ are $5\Z$ each.
Therefore, $\mO_{S(\Q)}=\Z[\frac{1}{5}]$ and hence 
${\rm rank}(\mO_{S(\Q)}^*)=1$.
This is an example of a set of primes of the CM-field $K$ which satisfies
the hypothesis.
\end{example}

We need Corollary \ref{goodunit2} to prove the main theorem of the paper.
\section{Proof of the Main Theorem}
We imitate the proof of the same
result for the case of arithmetic subgroups of $SL_2(K)$ 
(cf. \cite{mine2}). Here, we state a result due to Vaserstein which we
use in the proof of Theorem \ref{maintheorem}.

\begin{theorem} [Vaserstein] \cite{Vaserstein}\label{TheoremVaserstein}
Let $K$ be a number field and let $S$ be a finite set of primes in $K$
including $S_\infty$ such that $\card(S)\ge2$. Let $\mf{a}$ 
be a nonzero ideal of $\mO_S$. The group generated by
$\begin{pmatrix}1&\mf{a}\cr0& 1\end{pmatrix}$
and $\begin{pmatrix}1&0\cr \mf{a}&1\end{pmatrix}$
is a subgroup of finite index in $SL_2(\mO_S)$.
\end{theorem}

To prove Theorem \ref{maintheorem}, it suffices to show that 
any subgroup of finite index in $SL_2(\mO_S)$ is virtually generated by
three elements. 
Let $\Gamma$ be a subgroup of finite index in $SL_2(\mO_S)$. Without
loss of generality we assume that it is a normal subgroup. Let
its index in $SL_2(\mO_S)$ be $h$.\\

\noindent
{\bf Case 1:}
The pair $(K,S)$ is such that for {\bf every proper subfield} 
$L$ of $K$, we have
\begin{equation}\label{rankinequality}
{\rm rank}(\mO_{S(L)}^*\cap \mO_S^*)<{\rm rank}(\mO_S^*).
\end{equation}

Choose $\alpha\in \mO_S^*$ {\bf as in Corollary \ref{goodunit2}}. 
Obviously, $\begin{pmatrix}\alpha^h&0\cr0&\alpha^{-h}\end{pmatrix}\in \Gamma$.
Since $\Z[\alpha^h]$ is a subring of finite index in $\mO_S$, 
we replace $\alpha^h$ by $\alpha$ and assume that 
$\gamma:=\begin{pmatrix}\alpha&0\cr0&\alpha^{-1}\end{pmatrix}\in \Gamma$.
Define, $\psi_1:=\begin{pmatrix}1&0\cr h&1\end{pmatrix}\in \Gamma$ and
$\psi_2:=\begin{pmatrix}1& h\cr 0&1\end{pmatrix}\in \Gamma$. 
Let $\Gamma_0=\big<\gamma,\psi_1,\psi_2\big>$. We claim that 
$\Gamma_0$ is a subgroup of finite index in $SL_2(\mO_S)$.

Indeed,
$\gamma^{-r}\psi_1^s\gamma^r=\begin{pmatrix}1& 0\cr s\alpha^{2r}h &1
\end{pmatrix}\in \Gamma_0$ and $\gamma^r\psi_2^s\gamma^{-r}=
\begin{pmatrix}1& s\alpha^{2r}h \cr 0&1\end{pmatrix}\in \Gamma_0$. 
One concludes from this that $\Gamma$ contains 
$\begin{pmatrix}1&x\cr 0&1 \end{pmatrix}$ and
$\begin{pmatrix}1&0\cr y&1 \end{pmatrix}$ with $x,y\in h\Z[\alpha^2]$.
By Corollary \ref{goodunit2}, $h\Z[\alpha^2]$ is of 
finite index in the additive group $\mO_S$. If $m$ is the index then the ideal
$\mf{a}:=m\mO_S$ is contained in $h\Z[\alpha^2]$. Now it follows from Theorem 
\ref{TheoremVaserstein} that the group 
$\Gamma_0$ is a subgroup of finite index in 
$SL_2(\mO_S)$.\\

\noindent
{\bf Case 2:}
The pair $(K, S)$ is such that 
the inequality \eqref{rankinequality} {\bf does
not hold for some proper subfield} $F$ of $K$. That is, we have 
\begin{equation}\label{case2}
\rank(\mO_{S(F)}^*\cap\mO_S^*)=\rank(\mO_S^*). 
\end{equation}

This implies that $\rank(\mO_{F}^*)=\rank(\mO_K^*)$. Thus,  
by Lemma \ref{infiniteindex}, $K$ is a CM field and in fact 
$K=F(\sqrt{-d})$ so that $F$ is a totally real field
and $d$ a totally positive integer in $F$. 
Thus, we have
\begin{equation}\label{rankequality}
\mO_{S(F)}^*\asymp \mO_S^*,
\end{equation}
\begin{equation}\label{groupequality1}
\mO_F^*\asymp \mO_K^*.
\end{equation}
We prove a number theoretic lemma here. 
\begin{lemma}
With the above notation, let  \eqref{case2} hold for a 
CM filed $K=F[\sqrt{-d}]$. There exists $\alpha\in \mO_{S(F)}^*\cap\mO_S^*$ 
such that the ring $\Z[\alpha^n][\sqrt{-d}]$ is of finite index in $\mO_S$ 
for any integer $n$.
\end{lemma} 

\begin{proof}
In the case of a quadratic extension, a prime ideal 
of the base field is 
either inert or totally ramified or split completely (into two distinct 
primes). We claim
that the set $S(F)$ (cf. \eqref{definition1}),  does not contain
any finite prime which splits completely in $K$. 
To the contrary, if $S(F)$ contains a split 
prime $\mf{q}$ so that $\mf{q}\mO_K=\mf{q}_1\mf{q}_2$, then we have two 
possibilities, namely, $\mf{q}_1,\mf{q}_2\in S$ or $\mf{q}_1\in S$
and $\mf{q}_2\notin S$. If
$\mf{q}_1,\mf{q}_2\in S$, then $\card(S(F))<\card(S)$ (since $\mf{q}_1$
and $\mf{q}_2$ are contracted to the same prime $\mf{q}$ in $F$) and thus
\eqref{case2} does not hold. This is a contradiction. 
Next, assume that $\mf{q}_1\in S$ and $\mf{q}_2\notin S$. Let
$\beta$ (resp. $\gamma_1$) be the generator of $\mf{q}^{\ord_F(\mf{q})}$
(resp. $\mf{q}_1^{\ord_K(\mf{q}_1)}$). By \eqref{rankequality}, 
we have $\mO_S\supset\mO_{S(F)}$ so that 
$\beta\in\mO_S$. Again \eqref{rankequality} and 
\eqref{groupequality1} together imply
that $\gamma_1^m\in \mO_{S(F)}$ for some
$m>0$ so that $\gamma_1^{m}=u\beta^bx$ for some $b>0$ and
$u\in\mO_K^*\cap\mO_F^*$ and $x\in \mO_{S(F)}^*\cap\mO_S^*$ with 
$v_{\mf{p}}(x)=0$ for $\mf{p}\notin S(F)$.
Then $v_{\mf{q}_2}(\gamma_1)=0$ whereas
$v_{\mf{q}_2}(u\beta^bx)>0$. This is a contradiction.
Therefore, we have 
\begin{equation}\label{nonsplit}
(\mf{q}\cap\mO_F)\mO_K=\mf{q} \textnormal{ or }\mf{q}^2.  
\end{equation}
Let $\ord_F(\mf{q}\cap\mO_F)=a$. Then, by \eqref{nonsplit}, we see that
$(\mf{q}\cap\mO_F)^a\mO_K=((\mf{q}\cap\mO_F)\mO_K)^a=\mf{q}^a$ or 
$\mf{q}^{2a}$ is a principal ideal. Thus, $(\mf{q}\cap\mO_F)^a$ and
$\mf{q}^b$ (for $b=a$ or $2a$) are generated by the same element 
$\beta\in\mO_F$. 

Let 
$S-S_\infty=\{\mf{p}_1,\dots,\mf{p}_s\}$. 
Choose $\beta_i\in \mO_F$ such that 
$(\mf{p}_i\cap\mO_F)^{\ord_F({\mf{p}_i}\cap\mO_F)}=\beta_i\mO_F$. 
Then, by Lemma \ref{OS}, 
$\mO_{S(F)}=\mO_F[\beta_1^{-1},\dots,\beta_m^{-1}]$ and 
$\mO_S=\mO_K[\beta_1^{-1},\dots,\beta_m^{-1}]$.
Now, since $\mO_F[\sqrt{-d}]$ is of finite index in $\mO_K$, 
by Lemma \ref{comparison}, we have
$\mO_{S(F)}[\sqrt{-d}]=\mO_F[\sqrt{-d}][\beta_1^{-1},\dots,\beta_m^{-1}]$ 
is of finite index in $\mO_S$. Since $F$ is a non-CM field, by Theorem 
\ref{goodunit}, one can
choose $\alpha\in \mO_{S(F)}^*\cap\mO_S^*$ such that $\Z[\alpha^n]$ is of 
finite index in $\mO_{S(F)}$ for every $n\ge1$. Then 
$\Z[\alpha^n][\sqrt{-d}]$ is of finite index in $\mO_S$.
\end{proof}

Define $\gamma$ and $\psi_1$ 
as in {\bf case 1}. We modify the definition of $\psi_2$ by 
$\psi_2:=\begin{pmatrix}1& h\sqrt{-d}\cr0&1\end{pmatrix}\in \Gamma$. 
Let $\Gamma_0:=\big<\gamma,\psi_1,\psi_2\big>$. We show that $\Gamma_0$ 
is a subgroup
of finite index in $SL_2(\mO_S)$. 

Since $F$ is a non-CM field, by an
argument similar to case 1, one  
shows that there is an ideal $\mf{a}$ of $\mO_{S(F)}$ such that 
\begin{equation}\label{basicnonCM}
\begin{pmatrix}
1&0\cr\mf{a}&1\end{pmatrix}\subset \Gamma_0~\textnormal{ and }~
\begin{pmatrix} 1& \sqrt{-d}\mf{a}\cr0&1\end{pmatrix}\subset\Gamma_0.
\end{equation}
Then for $x\in \mf{a}$, using  Bruhat decomposition (see \cite[8.3]{Springer})
of $\psi_2$, we have 
\begin{equation}\label{bruhat1}
\begin{pmatrix}\hskip-20pt^{\psi_2}\hskip+10pt
1&0\cr x&1\end{pmatrix}
=~\begin{pmatrix}\cgtu\hskip-5pt 1&
h^2dx\cr0&1\end{pmatrix}\in\Gamma_0
\textnormal{ where }u=\begin{pmatrix}
 1&0\\\frac{1}{h\sqrt{-d}}&1\end{pmatrix}. 
\end{equation}

\noindent
Let $\mf{b}=h^2d\mf{a}$. Then we have
\begin{equation}\label{ucontainment}
\begin{pmatrix}\cgtu\hskip-5pt 1& \mf{b}
\cr 0&1\end{pmatrix}\subset\Gamma_0 ~~{ \rm and } ~~ 
\begin{pmatrix}\cgtu\hskip-5pt 1& 0\cr
\mf{b}&1\end{pmatrix}=\begin{pmatrix}1& 0\cr
\mf{b}&1\end{pmatrix}\subset\Gamma_0.
\end{equation}

\noindent
Let $\Gamma_1$ be the subgroup of $SL_2(\mO_F)$ generated by
$\begin{pmatrix}1& \mf{b}\cr0&1\end{pmatrix}$ and $\begin{pmatrix}1& 0
\cr\mf{b}&1\end{pmatrix}$. Then, by \eqref{ucontainment}, we have  
$^u\Gamma_1\subset\Gamma_0$. By Theorem \ref{TheoremVaserstein}, the index 
of $\Gamma_1$ in $SL_2(\mO_F)$ is finite. Thus it follows that 
there exists an integer $N$ 
such that 
\begin{equation}\label{belong}
\gamma^N\in\Gamma_1\cap\Gamma_0. 
\end{equation}
Since 
$^u\Gamma_1\subset\Gamma_0$, we have
$^u\gamma^N\in\Gamma_0$.

Therefore, 
$^u\gamma^{-N}\gamma^N=
\begin{pmatrix}1& 0 \cr (\alpha^{2N}-1)\frac{\sqrt{-d}}{hd}&1\end{pmatrix}\in
\Gamma_0$. Now by conjugating this element and its powers by 
negative powers of $\gamma$,
one shows that 
\begin{equation}\label{firstnonCM}
\Gamma_0\supset\begin{pmatrix}1&0\cr\sqrt{-d}\mf{c}&1
\end{pmatrix}
\end{equation}
where $\mf{c}:=(\alpha^{2N}-1)\Z[\alpha^2]\cap\mf{a}$.
Now $\mf{c}+\sqrt{-d}\mf{c}$ is a subgroup of finite index in 
$\mO_{S(F)}[\sqrt{-d}]$ and hence in $\mO_S$. 
Therefore, the group $\mf{c}+\sqrt{-d}\mf{c}$ 
contains a nonzero ideal $\mf{q}$ of $\mO_S$. 
Since $\mf{c}\subset\mf{a}$, by \eqref{basicnonCM} and 
\eqref{firstnonCM}, we have
\begin{equation}\label{firstcontainment}
\begin{pmatrix}1&0\cr\mf{q}&1
\end{pmatrix}\subset \Gamma_0. 
\end{equation}
\noindent
Again, for $y\in \mf{a}$, using the Bruhat decomposition of $\psi_1$, we have
\begin{equation}\label{bruhat2}
\begin{pmatrix}\hskip-20pt^{\psi_1}\hskip+10pt 1&y\sqrt{-d}\\0&1\end{pmatrix}
=\phantom{....}\begin{pmatrix}\hskip-27pt^{v\varphi}~\hskip+10pt
1&0\\h^2yd&1\end{pmatrix}\in \Gamma_0
\end{equation}
where $v=\begin{pmatrix}1&\frac{1}{h}\cr 0&1\end{pmatrix}$ and 
$\varphi=\begin{pmatrix}1&0\cr 0&\frac{1}{\sqrt{-d}}\end{pmatrix}$. 
Thus we have
\begin{equation}\label{uphicontainment}
~~~\begin{pmatrix}\hskip-27pt^{v\varphi}~\hskip+10pt1& \mf{b}
\\0&1\end{pmatrix}\subset\Gamma_0,\phantom{.........}\textnormal{and}
\phantom{.......}
\begin{pmatrix}\hskip-27pt^{v\varphi}~\hskip+10pt 1& 0\\
\mf{b}&1\end{pmatrix}\subset\Gamma_0.
\end{equation}
Therefore, $^{v\varphi}\Gamma_1\subset\Gamma_0$ and hence
$^{v\varphi}\gamma^N\in\Gamma_1\cap\Gamma_0$.  Thus, using \eqref{belong}
we have 
\begin{equation}
^{v\varphi}\gamma^N\gamma^{-N}=
\begin{pmatrix} 1&(1-\alpha^{2N})\frac{1}{h}
\\0& 1\end{pmatrix}\in\Gamma_0.
\end{equation}
Again by conjugating this element and its powers by nonnegative powers of 
$\gamma$, one shows that 
\begin{equation}\label{1containment}
\begin{pmatrix}1&\mf{c}\cr0&1\end{pmatrix}\subset\Gamma_0.
\end{equation}

\noindent
Since 
$\mf{c}\subset\mf{a}$, by \eqref{basicnonCM} and \eqref{1containment}, we have
\begin{equation}\label{secondcontainment}
\begin{pmatrix} 1& \mf{q}
\\0& 1\end{pmatrix}\subset\Gamma_0.
\end{equation}
It follows from \eqref{firstcontainment} and \eqref{secondcontainment}, and
by Theorem \ref{TheoremVaserstein}, that the group $\Gamma_0$ is a
subgroup of finite index in $SL_2(\mO_S)$. 
This completes the proof of Theorem \ref{maintheorem}.
\hfill$\square$


\end{document}